\documentclass[twoside]{article}
\usepackage{amsfonts}
\usepackage{amsmath}
\usepackage{amsthm}
\usepackage{amssymb}
\usepackage{amscd}
\usepackage{mathrsfs}
\usepackage{graphicx}
\usepackage{psfrag}
\usepackage{fancyhdr}
\usepackage{times}

\newcommand{\ra}{\rightarrow}
\newcommand{\hra}{\hookrightarrow}
\newcommand{\xra}{\xrightarrow}

\newcommand{\R}{\mathbb{R}}
\newcommand{\Q}{\mathbb{Q}}

\newcommand{\CP}{\mathbb{P}}

\newcommand{\Li}{\mathcal{L}_i}

\newcommand{\Tgn}{\mathcal{T}_{g,n}}

\newcommand{\M}{\mathcal{M}}
\newcommand{\Mb}{\mathcal{\overline{M}}}
\newcommand{\Mc}{\mathcal{M}^{comb}}
\newcommand{\Mgnc}{\mathcal{M}_{g,n}^{comb}}

\newcommand{\Mgnd}{\mathcal{M}_{g,n}^{\Delta}}
\newcommand{\Md}{\mathcal{M}^{\Delta}}
\newcommand{\Mcb}{\mathcal{\overline{M}}^{comb}}
\newcommand{\Mgncb}{\mathcal{\overline{M}}_{g,n}^{comb}}

\newcommand{\Mgn}{\mathcal{M}_{g,n}}
\newcommand{\Mgnb}{\mathcal{\overline{M}}_{g,n}}

\newcommand{\Mg}{\mathcal{M}_{g}}

\newcommand{\He}{\mathcal{H}}
\newcommand{\Hg}{\mathcal{H}_g}
\newcommand{\Hgb}{\mathcal{\overline{H}}_g}
\newcommand{\Hb}{\mathcal{\overline{H}}}

\newcommand{\Hgc}{\mathcal{H}^{comb}_g}

\newcommand{\GHg}{\Gamma_{H_g}}

\newcommand{\Wb}{\overline{W}}

\providecommand{\abs}[1]{\lvert#1\rvert}

\providecommand{\dof}[1]{\mathrm{d}#1\,}

\newcommand{\homeo}{\cong}
\newcommand{\iso}{\cong}
\newcommand{\define}{\equiv}

\newcommand{\Teich}{Teichm\"uller }
\newcommand{\Poin}{Poincar\'e }
\newcommand{\Weier}{Weierstra\ss{} }

\newcommand{\Aut}{\mathrm{Aut}}

\newcommand{\oSim}{\Delta^{\circ}}
\newcommand{\oOrb}{\oSim_{\Gamma}/\mathrm{Aut}(\Gamma)}

\newcommand{\varep}{\varepsilon}

\newcommand{\be}{\begin{equation}}
\newcommand{\ee}{\end{equation}}

\newtheorem{thm}{Theorem}[section]
\newtheorem*{mainthm}{Theorem} 
\newtheorem{lemma}[thm]{Lemma}
\newtheorem*{cor}{Corollary}

\theoremstyle{definition}
\newtheorem*{defn}{Definition}

\theoremstyle{remark}

\author {Alex James Bene}
\title {Combinatorial Classes, Hyperelliptic Loci, and Hodge Integrals}
\date{}

\begin{document}
\renewcommand{\headrulewidth}{0pt}
\renewcommand{\headsep}{15pt}
\maketitle

\begin{abstract}
A closed formula is obtained for the integral $\int_{\Hgb^1}\kappa_{1}\psi^{2g-2}$ of tautological classes over the locus of hyperelliptic \Weier points in the moduli space of curves.  As a corollary, a relation between Hodge integrals is obtained. 
  
The calculation utilizes the homeomorphism  between the moduli space of curves $\M_{g,1}$ and the combinatorial moduli space $\Mc_{g,1}$, a PL-orbifold whose cells are enumerated by fatgraphs.  This cell decomposition can be used to naturally construct combinatorial PL-cycles $W_a\subset\Mc_{g,1}$ whose homology classes are essentially the \Poin  duals of the Mumford-Morita-Miller classes $\kappa_a$.    In this paper we construct another  PL-cycle   $\Hgc \subset \Mc_{g,1}$ representing  the locus of hyperelliptic \Weier points and explicitly describe the chain level intersection of this cycle with $W_1$.   Using this description of $\Hgc\cap W_1$,  the  duality between  Witten cycles $W_a$ and  the $\kappa_a$ classes, and  Kontsevich's scheme of  integrating $\psi$ classes, the   integral $\int_{\Hgb^1}\kappa_{1}\psi^{2g-2}$ is reduced  to a weighted sum over graphs  and is evaluated by the enumeration of trees. 
\end{abstract}

\section {Introduction}
\label{sect:intro}

The moduli space of curves $\Mgn$ and its compactification $\Mgnb$ are objects of central interest, and the study of their tautological rings $R^\bullet(\Mgn)$ and $R^\bullet(\Mgnb)$ has been one of the primary avenues for probing their structure.  In particular, relations of top dimensional classes in $R^\bullet(\Mgnb)$, or intersection numbers, have been greatly explored in  the past two decades.
Although the tools of algebraic geometry have proven indispensable for a complete study of the moduli space of curves, it seems that most geometric results are eventually reduced to topological properties, and the use of topological methods has also contributed greatly to the field.

One of the most useful topological tools for probing the structure of the moduli of curves  has been its combinatorial description inspired by the insights of Thurston and Mumford. For the case of one marked point,  this description gives a cell decomposition of $\M_{g,1}$ in terms of arc families in a genus $g$ surface, which  Harer and Zagier first used  to determine the orbifold Euler characteristic of $\M_{g,1}$  \cite{harerzagier}. 

 More generally, the combinatorial description  gives a homeomorphism (in both the conformal \cite{harer} and hyperbolic \cite{penner1} settings) between $\Mgn\times\oSim_{n-1}$ and a combinatorially defined PL-orbifold,  the combinatorial moduli space $\Mgnc$.  Penner's subtle reinterpretation of this cell decomposition of $\Mgn\times\oSim_{n-1}$  and its combinatorics in terms of Feyman-like fatgraphs \cite{penner2} gave a new perspective on Harer and Zagier's result and lead to a more natural 
connection with the matrix models of physics.  In particular, the fatgraph-enumerated cells  transform  integrals over the moduli space $\Mgn$ (or $\Mgnb$)  into Feynman-type  weighted sums of fatgraphs.

Perhaps the most important culmination of the fatgraph-matrix model perspective has been Kontsevich's proof of Witten's conjecture which provided a recursive method of determining all intersection products of descendant $\psi$ classes in $R^\bullet(\Mgnb)$ based on the equations of the KdV hierarchy \cite{kontsevich}.  Kontsevich's proof relied heavily on the PL-orbifold structure and cell decomposition of $\Mgnc$, and much in the spirit of  rational homotopy theory, his approach made critical use of PL-forms with respect to this decomposition  to capture the relevant ring structure of $R^\bullet(\Mgnb)$ displayed by the $\psi$ classes.   In particular, by extending the fatgraph description to a (non-orbifold) compactification $\Mgn'$ of $\Mgn$, he employed PL-de Rham style integration techniques and the combinatorics of matrix models to determine intersection numbers of $\Mgnb$.  

\subsection{Results and Methods}
\label{sect:results}
Although Kontsevich's result in principle determined all intersection numbers of $\psi$ classes, closed formulae for integrals of tautological classes are in general not immediately derived.  In particular, the discovery of  closed formulae for top intersections  of $\psi$ and $\lambda$ classes, called Hodge integrals, is still an active field.  

 The main theorem of this paper  is the derivation of the following integral of tautological classes over the locus of hyperelliptic \Weier points
\begin{mainthm}
For $g\geq 1$,
\be\label{eq:intromain}
\int_{\Hgb^1}\kappa_1\psi^{2g-2}    =  \frac{(2g-1)^2}{2^{2g}(2g+1)!}
\ee
\end{mainthm}
which produces as a corollary the following relation of Hodge integrals
\begin{cor}
For $g\geq2$,
\be\label{eq:introcor}
\sum_{i=0}^{g-1} \int_{\Mb_{g,1}} \left[ (-1)^{i}(2^{g-i}-1)\psi^{g-i-1}\lambda_{i} \right]\kappa_1\psi^{2g-2} =  \frac{14g^2 -11g +3  }{3\cdot{}2^{2g}(2g+1)!}.
\ee
\end{cor}

Although the proof of the theorem relies on the enumerative combinatorics  of planar trees rather than matrix models, the methods used are  much in the spirit of the Witten-Kontsevich model and are essentially topological in nature.  In fact, the calculation can be performed, but not immediately justified, in terms of Kontsevich's  compactification $\M_{g,1}'$ alone.  The corollary, on the other hand, requires a more careful analysis of the Deligne-Mumford compactification $\Mb_{g,1}$ and a delicate application of Porteus' formula.

Perhaps the most striking difference between our approach  and  Kontsevich's  is that ours utilizes the intersection theory of $\Mb_{g,1}$ in terms of PL-cycles rather than exclusively in terms of PL-forms.  Thus instead of relying entirely on integration to determine the intersection numbers, our calculation makes more explicit use of the (virtual) \Poin duality of $\Mgnb$ and involves analyzing chain level intersections much in the vein of classical analysis situs.   As we will see, these chain level intersections are closely related to the facet structure of Stasheff's associahedron and can be thought of as a new source of combinatorics in the problem.

The key idea which allows us to reinterpret the integral \eqref{eq:intromain} in terms of cycles comes from the suggestion of Witten that the Mumford-Morita-Miller classes $\kappa_a$ might be (in some sense) \Poin dual to naturally defined combinatorial PL-cycles $W_a$ in $\Mgnc$.  These Witten cycles  $W_a$ are defined as certain collections of cells of (real) codimension at least $2a$ corresponding to fatgraphs with one vertex of valence at least $2a+3$.  Their duality with the $\kappa_a$ classes has been proven to various degrees of generality in various contexts.  For our result, we shall make essential use of the  $a=1$ case of the  extension of Witten's conjecture to  $\Mb_{g,1}$ which states  \cite{mondello} that under \Poin duality
\be \label{eq:introduality}
\kappa_1\cong \frac{1}{12}\left( [\Wb_1] + [\partial\Mgn] \right).
\ee

By applying the above duality \eqref{eq:introduality}, we can rewrite the integral \eqref{eq:intromain} of the theorem as a sum of two integrals
\[
\int_{\Hgb^1}\kappa_1\psi^{2g-2}    = \frac{1}{12}\left[ \int_{[\Hgb^1]\cdot [\Wb_1]}\psi^{2g-2} +\int_{[\Hgb^1]\cdot [\partial\M_{g,1}]} \psi^{2g-2} \right].
\]
The evaluation of the first integral on the RHS in the above expression is the heart of our calculation.  It involves explicitly describing the locus of the intersection  $\Hg^1\cap W_1$, which includes both
determining  the multiplicities of its components and enumerating its cells.  

To fully analyze the intersection $\Hg^1\cap W_1$, we first construct a PL-cycle $\Hgc\subset\Mc_{g,1}$  and argue that it corresponds to $\Hg^1$ by recognizing the symmetry of the hyperelliptic involution in certain fatgraphs.  One advantage of our description of $\Hg^1$ is that  its cells are easily  enumerated by certain decorated planar trees.  Thus, we are able to reduced the integral to a Feynman-type  sum over decorated trees and consequentially  evaluate it by the use of generalized Catalan numbers. 

\subsection{Comments}
\label{sect:acknowledge}
Part of our interest in the calculation of this paper is that it is in some sense the first application of the  duality suggested by Witten  to determining relations in the tautological ring.  More generally,  Witten's duality  can in principle translate  the problem of finding intersection numbers of $\kappa_a$  classes into questions about the intersections of the PL-cycles $W_a$.  However, as the intersections of Witten cycles $W_a$ are highly degenerate, answering these questions would require  some form of chain level general positioning argument and appears to be difficult.  It is our hope that the methods of this paper  might give some insight into and serve as a first step towards this more difficult problem of general combinatorial intersections.

We also feel that the PL description $\Hgc$ of the locus of hyperelliptic \Weier points given in this paper may be of interest in its own right.  In particular, it  may be useful for exploring properties of the hyperelliptic mapping class group.  

While we feel that the methods of this paper are important for a full understanding of the results, especially their connections with matrix models, we briefly  mention here two other approaches that are relevant.   Firstly, the theorem of this paper has also been  derived and generalized by J. Bertin and M. Romagny using their more algebro-geometric treatment of the locus of hyperelliptic \Weier points \cite{bertin}.  
Secondly, it appears that the use of single Hurwitz numbers  should lead to another proof of the corollary. 
   We hope to explore the connections between our techniques and these in future papers.  

 


\section{The Moduli Space of Curves}
\label{sect:moduliofcurves}

The moduli space of $n$ pointed genus $g$ curves $\Mgn$ as a set is the set  equivalence classes $\{(C,x_1,\ldots,x_n)\}/\!\sim$ of genus $g$ compact Riemann surfaces, or smooth algebraic curves,  marked by $n$ distinct ordered points, where the equivalence is given by biholomorphic maps taking ordered marked points to ordered marked points.  We shall always assume the condition $2g-2+n>0$, in which case a natural topology is given to $\Mgn$ by its construction as the quotient of \Teich space $\Tgn$, a contractible $6g-6+2n$ (real) dimensional space, by the properly discontinuous action of the mapping class  group $Mod_{g,n}$.  In this way, $\Mgn$ has the structure of a differentiable orbifold with rational cohomology  isomorphic to that of $Mod_{g,n}$.

Algebro-geometrically, the moduli of curves is a quasi-projective variety with only quotient singularities but is more often referred to by algebraic geometers as a smooth Deligne-Mumford stack.  For topologists, the greatest  advantage of the algebro-geometric viewpoint is that it leads to a natural orbifold compactification $\Mgnb$, the moduli of stable $n$-pointed genus $g$ algebraic curves.  Here,  an $n$-pointed genus $g$ stable curve is an algebraic curve of arithmetic genus $g$ with only double point singularities together with $n$ ordered  smooth marked points such that each component $C_i$ satisfies the stability condition: $2g_i -2 + n_i>0$, where $g_i$ is the genus of the component $C_i$ and $n_i$ denotes number of marked  points plus the number of nodes of  $C_i$.  

The importance of the orbifold compactification $\Mgnb$ is that it can be studied using many of the tools familiar from compact manifolds.  In particular, $\Mgnb$ satisfies (virtual) \Poin duality, thus the ring structure of its (de Rham) cohomology ring translates into a proper intersection theory in terms of topological cycles. We shall not give here the formal definition of orbifolds beyond the fact that they are locally modeled on quotients of a Euclidean ball under the action of a finite group, and we will rather naively treat them more or less as manifolds with two caveats.  Firstly, the fundamental class of an orbifold cycle will always be weighted by dividing by the order of its generic symmetry group.  Secondly, the multiplicity of  intersecting orbi-cycles will be calculated by locally lifting the cycles to their Euclidean neighborhoods, where the intersection can be viewed classically.  From this point forward we shall almost always drop the prefix `orbi' when referring to simplices, cells, cyles, etc.  

\subsection{The Tautological Ring}

We now  review the basic definitions of the tautological classes and tautological ring of $\Mgnb$.  For a more detailed exposition, we refer to a few of the many excellent resources \cite{arbarello}\cite{kontsevich}\cite{harris}\cite{witten}\cite{mondello}.  We note here that although the definition of the tautological classes and tautological ring are made in the context of the Chow rings $A^\bullet(\Mgnb,\Q)$ of algebraic geometry, we shall be primarily interested in their images in the rational cohomology under the cycle class map.  From this point forward, all coefficients are assumed taken in $\Q$.  

We begin by describing several natural proper morphisms associated with the moduli space of curves.  First, there is the  flat \emph{forgetting} morphism
\be\label{eq:forgetting}
 \pi_{n+1}: \Mb_{g,n+1} \ra \Mgnb
 \ee
  which over a smooth curve simply forgets the $n+1$st marked point and over a singular curve forgets the point and contracts any non-stable components (which are necessarily twice marked spheres).  The map $\pi_{n+1}$  identifies $\Mb_{g,n+1}$ as the ``universal curve'' over $\Mgnb$, and there exist  $n$ natural sections, or \emph{bubbling up} morphisms $s_i: \Mgnb\ra\Mb_{g,n+1}$ obtained by adding the $n+1$st marked point over the $i$th marked point, thus resulting in a rational component marked by a node, $x_i$, and $x_{n+1}$.   If we let $\omega_{\pi_{n+1}}$ denote the relative dualizing sheaf  of $\pi_{n+1}$ on $\Mgnb$  (an extension of the sheaf of sections of the cotangent bundle along the fibers over  $\Mgn$), and denote by $D_i$ the divisors defined by the images of the $s_i$,  we can then define the tautological classes as follows:
\begin{equation}\label{tautdef}
\begin{split}
\psi_i  & =    c_1(s_i^*\omega_{\pi_{n+1}}) \\
\kappa_a  & =  \pi_{n+1*}(c_1\big(\omega_{\pi_{n+1}}\big(\sum{}D_i\big)\big)^{a+1}) \\
\lambda_j & =  c_j(\pi_{n+1*}(\omega_{\pi_{n+1}})) 
\end{split}
\end{equation}
where $c_j$ denotes the $j$th Chern class, and $\pi_{n+1*}$ denotes the push forward morphism associated with $\pi_{n+1}$, which in the topological setting is just the Umker map.  

The geometric meanings of the tautological classes defined by \eqref{tautdef} can be described through naturally defined bundles over $\Mgnb$ as follows.
\[
\begin{CD}
\mathcal{L}_i 	@<<<	T^*_{x_i}C  @. \quad  @.    \quad    @.	   \mathbb{E}  @<<<   	H^0(C,\Omega) \\
@VVV  		@.  				@.   @.	@VVV	@.		@.   \\
\Mgnb		@.		@.		@.   @. 	\Mgnb	@.		@.
\end{CD}
\]
Let $\Li$ be the line bundle over $\Mgnb$ such that the fiber over a point $(C, x_1,\dotsc, x_n)\in\Mgnb$ is the cotangent space $T^*_{x_i}C$ of $C$ at $x_i$.  The descendent classes $\psi_i\in A^1(\Mgnb)$  (resp.\ $\in H^2(\Mgnb)$) are then the Chern classes of these bundles, $\psi_i=c_1(\Li)$.  
The Mumford-Morita-Miller classes $\kappa_a\in A^a(\Mgnb)$ (resp.\ $\in H^{2a}(\Mgnb)$),  are then simply the pushforwards, or integration along the fiber,  of powers of $\psi$ classes, 
\be \label{eq:kappadef}
	\kappa_a=\pi_{n+1*}(\psi_{n+1}^{a+1}).
\ee

Now, define the Hodge bundle $\mathbb{E}$ as the vector bundle over $\Mgnb$ with fiber over a point $(C, x_1,\dotsc, x_n)\in\Mgnb$ equal to $H^0(C,\omega_C)$, the sections of the dualizing sheaf (which for a  smooth curve is just  the vector space of holomorphic 1-forms).  The $\lambda$ classes are then the Chern classes of this bundle $\lambda_j=c_j(\mathbb{E})$.  

There also exist natural proper \emph{gluing} morphisms 
\be \label{eq:boundarymorphisms}
\delta^n_g \colon \Mb_{g,n+2} \rightarrow \Mb_{g+1,n} \quad \delta^{n,n'}_{g,g'}\colon  \Mb_{g,n+1} \times \Mb_{g',n'+1} \rightarrow \Mb_{g+g',n+n'}
\ee
where $\delta^n_g$ is the map which identifies the last two marked points of a genus g curve into a single  node and $\delta^{n,n'}_{g,g'}$ is the map which glues two curves of genera $g$ and $g'$ together along their last marked points.  Classes  which are pushforwards under the gluing morphisms \eqref{eq:boundarymorphisms} are called boundary classes. The total boundary  $\partial\Mgn$ itself represents a complex codimension one boundary class which will be particularly important for our calculation.

The tautological rings $R^\bullet(\Mgnb)$ are defined as the minimal system of $\Q$-algebras containing the tautological classes which are closed under the forgetting and boundary morphisms \eqref{eq:boundarymorphisms} and \eqref{eq:forgetting}. Again, although initially defined as a subring of the Chow ring, we shall primarily consider the tautological ring a subring of the cohomology ring by use of the cycle class map.  The tautological ring of open moduli $R^\bullet(\Mgn)$ is defined to be the restriction of $R^\bullet(\Mgnb)$ to $\Mgn$, and we shall not distinguish in notation the differences between the tautological classes $\psi_i$, $\kappa_a$, and $\lambda_j$ in these different rings.

As already mentioned, understanding the structure of the tautological ring is of considerable importance, and much of the work on the tautological rings has focused on relations of top dimensional  classes.   One well known relation between tautological classes that we shall make reference to  is  the string equation given by (see \cite{witten})
\be\label{eq:string}
\pi_{n+1*}(\prod_{i=1}^n \psi_i^{a_i}) = \sum_{j=1}^n \prod_{i=1}^n \psi_i^{\textrm{max}\{a_i-\delta_{ij},0\}}.
\ee
 One can check that this relation immediately leads to  the determination of all genus zero intersection numbers of $\psi$ classes and in particular gives 
\[
\int_{\Mb_{0,n}}\psi^{n-3}=\int_{\Mb_{0,3}}\psi^{0}=1
\]
as $\Mb_{0,3}$ is a single point without automorphism.   

\section{The Combinatorial Moduli Space}
\label{sect:mgnc}

As mentioned in the introduction, 
one of the most useful topological tools for probing the structure of the moduli of curves for $n>0$ is its description in terms of fatgraphs, given in the hyperbolic setting by Penner's decorated \Teich space \cite{penner1} and in the conformal setting by the Harer-Mumford-Strebel model \cite{harer} \cite{kontsevich}. 

\subsection{Fatgraphs}
\label{sect:fatgraphs}
A fatgraph is a finite connected graph, i.e.\ finite connected  1-dimensional CW complex, with an additional \emph{fat} structure given by a cyclic ordering of the half-edges, or stubs, incident upon each vertex.  For  convenience, we shall always assume that all fatgraphs have vertices of valence at least three, unless otherwise stated.  

\begin{figure}[!h]
\begin{center}
\psfrag{1}{$e_1$}
\psfrag{2}{$e_2$}
\psfrag{3}{$e_3$}
\psfrag{4}{$e_4$}
\psfrag{5}{$e_5$}
\includegraphics[width=4.4 in]{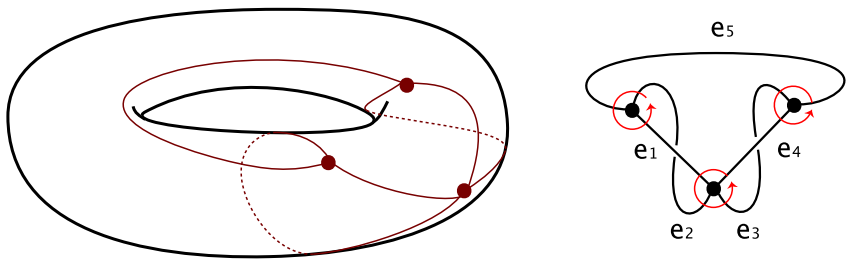}
\caption{Fatgraph of type $(1,2)$ minimally embedded in a surface.}
\label{embeddedfatgraph}
\end{center}
\end{figure}

Given  a fatgraph, one can consider embeddings  of the graph into a closed (or punctured) oriented surface, where the cyclic ordering of vertices are forced to coincide with the orientation of the surface. In this way, a fatgraph $\Gamma$ is often depicted as a ribbon graph which  is a cellular thickening of its  image in the surface.  To any fatgraph $\Gamma$, we can associate two numbers $g(\Gamma)$ and $n(\Gamma)$ where $g(\Gamma)$ is the minimal genus of surface in which $\Gamma$ can be embedded, and $n(\Gamma)$ is the number of contractible 2-cells in the compliment of $\Gamma$ under such an embedding.  We call a fatgraph of type $(g,n)$ if $g(\Gamma)=g$ and $n(\Gamma)=n$.

Given $\Gamma$ and any minimal genus embedding $\Gamma\hra S$, if we label the $n$ complimentary 2-cells $b_i$ for $i=1,\ldots,n$, then to each  $b_i$, we can associate the cyclicly ordered sequence of edges of $\Gamma$ which bound the  cell.  We  call this sequence the $i$th boundary cycle of $\Gamma$.   The fatgraph in figure \ref{embeddedfatgraph} has two boundary cycles $(e_1,e_4,e_3,e_2)$ and $(e_1,e_5,e_4,e_3,e_5,e_2)$. In all future discussions, we shall always implicitly assume that a fatgraph comes equipped with  a labeling of its boundary cycles, but when we wish to be explicit, we shall refer to such a pair as a boundary-labeled fatgraph.

We now discuss the types of maps that exist between fatgraphs.  Firstly, 
an isomorphism of (boundary-labeled) fatgraphs is a bijection of vertices and directed edges which preserves the incidence relations, the fat structure, and the labeling of the boundaries.  An automorphism is then just an isomorphism from a fatgraph to itself, which must by necessity be of finite order.  One can check that the fatgraph in figure \ref{embeddedfatgraph} has a single non-trivial automorphism of order two.

Given a non-loop edge $e$ in a fatgraph $\Gamma$, there is a map called the \emph{edge collapse} of $e$ to a new fatgraph $\Gamma_e$ with one less edge and one less vertex.  Specifically, we obtain $\Gamma_e$ by collapsing $e$, i.e.\  by  coalescing its endpoints $v_1$ and $v_2$ into a new vertex $v$  and giving the adjoining stubs of $v$   the naturally inherited cyclic ordering coming from those of $v_1$ and $v_2$. Note that the labeling of the boundaries is naturally preserved by such an edge collapse. In this situation, we also call $\Gamma$ an expansion of $\Gamma_e$.  If two (generally non-isomorphic) fatgraphs $\Gamma$ and $\Gamma'$ have edge collapses which are isomorphic, $\Gamma_e \iso {\Gamma'}_{e'}$,  we then say that they are related by a Whitehead move on $e$ (collapse of $e$ followed by expansion of $e'$). 

More generally, there exists maps which are compositions of isomorphisms and edge collapses.  The isomorphism type of the resulting fatgraph is uniquely determined by a collection of  edges to be collapsed, which consequentially must be a forest of disjoint trees. 

\subsection{The Combinatorial Moduli Space}
A metric fatgraph is a fatgraph $\Gamma$ together with an assignment of  a non-negative length $l(e)$  to each edge $e\in \Gamma$ such that no cycle in the fatgraph has total length zero.  Given a metric fatgraph $(\Gamma, l)$,  the total length of all edges in a boundary cycle $b_i$ is called the $i$th perimeter, and is denoted $p_i$.  The totality of all perimeters is denoted $p=(p_1,\ldots,p_n)$.    For a fixed  fatgraph $\Gamma$, the space of all   positive metrics on $\Gamma$ is homeomorphic to an open  cell 
 and the perimeter function $\rho$ maps this cell onto the interior of a cone in $\R_+^n$.  We call a metric \emph{normalized} if the sum of the perimeters equals one, i.e.\  the sum of lengths of all edges $\sum_e l(e)$ equals one half.   The space of all positive normalized metrics for a fixed $\Gamma$  can thus be identified with the interior of a simplex (after re-scaling by a factor of two), which we denote by $\Delta^\circ_\Gamma$.

Given an arbitrary metric on a fatgraph $\Gamma$, one can obtain a positive metric on a possibly different fatgraph by collapsing all edges of length zero.  This correspondence then identifies the interiors of certain faces of the simplex $\Delta^\circ_\Gamma$ with open simplices $\Delta^\circ_{\Gamma_i}$ of other fatgraphs.  Thus it  provides a method of gluing together simplices corresponding to different fatgraphs.  For convenience, when we refer to the underlying fatgraph $\Gamma$ of a metric fatgraph $(\Gamma,l)$, we shall always implicitly assume that the metric $l$ is positive on $\Gamma$. 

By an automorphism of a metric fatgraph $(\Gamma,l)$, we shall mean  the subgroup of the underlying fatgraph automorphism group $\Aut(\Gamma)$ which preserves the metric $l$. Obviously, the automorphism group of a fatgraph $\Gamma$ acts on its corresponding simplex $\oSim_\Gamma$ and we call the quotient under this action $\oOrb$ the orbi-simplex associated to $\Gamma$.  

We can now make the following
\begin{defn}
For  $g\geq 0$ and $n>0$ with $2g-2+n>0$, the combinatorial moduli space $\Mgnc$ of type $(g,n)$  is defined as  
\[
\Mgnc\define\coprod_\Gamma \oOrb
\]
where $\Gamma$ varies over all equivalence classes of boundary-labeled fatgraphs of type $(g,n)$, and the topology of $\Mgnc$ is given by the gluing described above.
\end{defn}

Note that the maximal cells of $\Mgnc$ are in 1-to-1 correspondence with  (equivalence classes of) trivalent fatgraphs of type $(g,n)$, thus a simple Euler characteristic argument shows $\Mgnc$ is a $6g-7+3n$ dimensional PL-orbispace (which is in fact an orbifold).  The perimeter function extends to $\Mgnc$ 
\[
\rho: \Mgnc \ra \oSim_{n-1}
\]
with image an open $n-1$-simplex.   

The main result regarding the combinatorial moduli space  is the following theorem which has essentially two proofs, one in the setting of hyperbolic geometry of punctured surfaces given by Penner's decorated \Teich space \cite{penner1}, and the other in the conformal setting based on the results of Strebel on quadratic differentials for Riemann surfaces \cite{strebel}\cite{harer}.
\begin{thm}[Penner, Strebel]  \label{thm:mainthm}
There exists a homeomorphism of orbifolds
\[
\Psi: \Mgnc \xra{\homeo} \Mgn\times\oSim_{n-1}
\]
with the second factor equal to the perimeter map.
\end{thm}

As a result, if we choose some $p\in\oSim_{n-1}$,  the \emph{slice}  $\Mgnc(p)\define\rho^{-1}(p)$  determined by $p$ gives a cell decomposition of $\Mgn$ itself
\[
\Mgn\homeo\coprod_\Gamma \oSim_\Gamma(p)/\Aut(\Gamma).
\]  
In this paper we shall primarily be interested in the case where $n=1$, where there is then a unique normalized slice, and we have 
\[
\Psi: \Mc_{g,1}\xra{\homeo} \M_{g,1}.
\]
Although the results of this paper will not depend on the explicit form of the homeomorphism $\Psi$,  from this point forward  we shall for concreteness assume that $\Psi$ is the map given by Strebel differentials.  


\subsection{Partial Compactifications}
It will be convenient to recall two partial compactifications of the combinatorial moduli space.  Both are obtained by allowing certain cycles of edges to collapse, much as we previously allowed non-loop edges (or more generally, trees of such edges)  to collapse.    The two compactifications are somewhat complimentary in the sense that the one allows for only boundary cycles to collapse, while the other allows only non-boundary cycles to collapse.  

The first compactification we consider is denoted $\Mgnd$, which is the partial compactification in the perimeter component.  The space $\Mgnd$ is built from cells of metric $\Delta$-labeled fatgraphs\footnote{We use the terminology $\Delta$-labeled rather than the established terminology $P$-labeled found in \cite{looijenga} and \cite{mondello}. The name was chosen to emphasize that  features  of the graph are  put into correspondence with the vertices of the perimeter simplex $\Delta_{n-1}$. 
}, 
an extension of boundary-labeled fatgraphs which  we  now define.  A $\Delta$-labeled fatgraph of type $(g,n)$  is a  fatgraph of type $(g,n')$ of arbitrary valence  with $n'<n$ labeled boundary cycles and $n-n'$ labeled vertices such that every vertex of valence less than three must be $\Delta$-labeled.  These labeled fatgraphs are naturally the limits of graphs where (all but one) boundaries are allowed to shrink to zero perimeter.  We refer to  \cite{looijenga} and \cite{mondello} for more details.  One advantage of this version is that there exists a natural slice given by the perimeter vector $(1,0,\ldots,0)\in\Delta_{n-1}$.

With this definition, we then have a generalization of theorem \ref{thm:mainthm},
\begin{thm}[Strebel]
There exists a homeomorphism of orbifolds
\[
\Psi^\Delta : \Mgnd \xra{\homeo} \Mgn\times\Delta_{n-1}
\]
extending the homeomorphism $\Psi$.
\end{thm}

The other partial compactification $\Mgncb$  of $\Mgnc$ is due to Kontsevich and includes cells corresponding to stable fatgraphs (for a definition, see \cite{kontsevich}), where cycles not containing any boundaries are allowed to collapse, thus resulting in possibly disconnected fatgraphs with special vertices labeled as nodes.  The corresponding theorem for this space is given by 
\begin{thm}[Kontsevich, Strebel]
There exists a homeomorphism of orbispaces
\[
\overline{\Psi} : \Mgncb \xra{\homeo} \Mgn'\times\oSim_{n-1}
\]
extending the homeomorphism $\Psi$, and 
where $\Mgn'$ is the quotient of $\Mgnb$ by the closure of the relation which identifies stable curves which are topologically equivalent but may differ in conformal structure on components which posses no marked points (see \cite{kontsevich} or \cite{zvonkine} for more details). 
\end{thm}

The spaces $\Mgncb$ and $\Mgn'$ are not orbifolds, but rather orbispaces which nonetheless posses well-defined fundamental classes.   Thus, top dimensional forms can be evaluated on these spaces; however, \Poin duality does not hold.  For the hyperbolic perspective on partial compactifications, we refer to \cite{penner5}.

\subsection{Expansions and Catalan Numbers}
For the main calculation of this paper, it will be essential that we understand the cell structure of a neighborhood of a point $(\Gamma,l)$ in $\Mc_{g,1}$.  For simplicity, we shall initially assume that $(\Gamma,l)$ has no automorphisms, thus the neighborhood is homeomorphic to a Euclidean ball.

  For  $(\Gamma,l)$ trivalent, a small neighborhood in $\Mgnc$ has a single cell comprised of metrics close to $l$ on the same fatgraph $\Gamma$.  However, for non-trivalent fatgraphs $(\Gamma,l)$, the neighborhoods consist of metrics on different expansions of $\Gamma$ (with the lengths of all expanded edges small).  Thus the neighborhood naturally decomposes into different cells, one for each expansion of $\Gamma$. 

The simplest case is the neighborhood of a metric fatgraph $(\Gamma,l)$ which is trivalent  except for a single $4$-valent vertex.   Every such $\Gamma$ corresponds to a particular Whitehead move and the two expansions of $\Gamma$ correspond to the two maximal cells which are glued together along their common  face $\oSim_{\Gamma}$.    

More generally, we shall be interested in the neighborhood of a metric fatgraph which is trivalent except for a single vertex of valence $k$ for $k>3$.  A small  neighborhood of such a metric fatgraph has cells corresponding to the ways of expanding the $k$-valent vertex and are in 1-to-1 correspondence with the equivalence classes  of rooted planar trees having $k$ leaves, or equivalently, the number of polygonizations of a regular $k$-gon. Thus the neighborhood has the cell structure of a Euclidean ball cross 
the \Poin dual fan structure of the associahedron $K_{k-1}$. 

The number of maximal cells in this neighborhood is equivalent to the number of rooted trivalent planar trees with $k-2$ non-leaf vertices (or triangulations of a $k$-gon) which is well known to be equal to the $k-2$nd Catalan number 
\be \label{eq:catalan}
C_{k-2}=\frac{1}{k-1}\binom{2k-4}{k-2}.
\ee
Perhaps the most direct  proof of this result can be obtained by using the branch decomposition of trivalent planar trees to obtain a relation on the generating function of the Catalan numbers. 

We shall also be interested in the number of partial expansions of a $k$-valent vertex which results in a graph with a single $5$-valent vertex, the rest trivalent (a total of  $k-4$ non-leaf vertices).  The number of these cells  is given by a generalization of the above Catalan numbers 
\be  \label{eq:gencatalan}
C_{5,k-4}=\binom{2k-6}{k-5}
\ee
as one can check either directly using \eqref{eq:catalan} and the branch decomposition of trees at their $5$-valent vertex, or by appealing to more general counting theorems (see \cite{goulden}).

A particular case which will be analyzed more fully throughout the paper is the neighborhood of a metric fatgraph that is trivalent except for a single $6$-valent vertex.  In this case, one can check that the corresponding dual associahedron has 14 3-cells corresponding to maximal expansions, 21 2-cells, and 9 1-cells, 6 of which preserve a $5$-valent vertex.

If the automorphism group of a metric fatgraph  $(\Gamma,l)$ is not trivial, 
the problem  of counting cells corresponding to fatgraph expansions is somewhat more subtle and difficult than the above counting arguments as  isomorphisms may exist between different expansions.  However, if one counts these cells in the orbifold sense by weighing each cell by the inverse of the order of its symmetry group, the answer  follows much as above.  In particular, for a metric fatgraph $(\Gamma,l)$  with a single $k$-valent vertex, we have the weighted sum of maximal cells $D_{\Gamma'}$ (in a small neighborhood of $(\Gamma,l)$) equal to 
\[
\sum_{\Gamma'>\Gamma} \frac{1}{\abs{\Aut(D_{\Gamma'})}} = \frac{C_{k-2}}{\abs{\Aut(\Gamma,l)}} 
\]
where the sum is over isomorphism classes $\Gamma'$ of trivalent expansions  of $\Gamma$ and $\abs{\Aut(\Gamma,l)}$ by necessity must divide $k$. 

\section{Integration and Combinatorial Classes}
\label{sect:integration}

One benefit of the combinatorial description of $\Mgn$ is that it provides a method of integration over $\Mgnb$.  In particular, given a top dimensional 
 form $\eta$ on $\Mgnb$ that has a description in terms of the fatgraph cells, we can evaluate the integral
\be \label{eq:intscheme}
\int_{\Mgnb} [\eta]=\int_{\Mgn} \eta=\sum_{\Gamma} \frac{1}{\abs{\Aut(\Gamma)}} \int_{\oSim_\Gamma(p)} \Psi^*(\eta)
\ee
as a weighted sum over trivalent fatgraphs by pulling back the form $\eta$ to each maximal cell ${\oSim_\Gamma(p)}$ of the slice $\Mgnc(p)$ for  some perimeter $p$.\footnote{In the case that $\eta$ is a differential form, the formula holds since the map $\Psi$ is smooth on the interior of cells.   Also, the fact that maximal cells of a slice correspond to slices of maximal cells, while non-obvious, has been assumed here.} Note that the factor $ \frac{1}{\abs{\Aut(\Gamma)}}$ is intuitive in the case where $\Aut(\Gamma)$ acts faithfully on the simplex slice $\oSim_\Gamma(p)$ since then the fundamental domain to be integrated over $\oSim_\Gamma(p)/\Aut(\Gamma)$ will be only that fraction of the whole simplex.  However, the above formula is valid even in the case when $\Aut(\Gamma)$ does not act faithfully, as it then contributes the proper symmetry factor for the orbifold fundamental class. 

This integration scheme \eqref{eq:intscheme} was first employed by Penner to compute Weil-Petersson volumes of $\Mgnb$ using a combinatorial description of the Weil-Petersson 2-form \cite{penner3}.  Kontsevich later made critical use of this technique in his proof of Witten's conjecture by finding a combinatorial PL representative of the $\psi$ classes which we now describe.

\subsection{Combinatorial  PL-Forms}
Let $b_k=(e_1,\ldots, e_m)$ be a boundary cycle of a possibly stable fatgraph $\Gamma$ (where some edges may repeat) and let $l_i=l(e_i)$ denote the simplicial coordinate corresponding to the metric on the edge $e_i$.  Let $\omega_k$ be the closed PL-form on $\oSim_\Gamma$ (or $\oOrb$) defined by
\[
\omega_k\define \sum_{1\leq i < j <m} \dof{\left(\frac{l_i}{p_k}\right)}\wedge\dof{\left(\frac{l_j}{p_k}\right)}.
\]

For a given perimeter value $p\in\Delta_{n-1}$, Kontsevich showed that the form $\Omega=\sum_i p_i^2 \omega_i$ restricted to the slice $\Mgnc(p)$ is symplectic and defines a compatible orientation on the maximal cells of the slice.  Moreover, for a fatgraph $\Gamma$ with only odd valence vertices (in particular, the trivalent ones corresponding to maximal cells), the volume form determined by $\Omega$ is related to that of the standard volume form on the cell $\oSim_\Gamma \times \R_+$ by 
  \[ 
 4^d2^{1-g}\prod_{e\in\Gamma}\abs{\dof{l(e)}}  = \frac{\abs{\Omega^{d}}}{d!}\prod_{i=1}^n \abs{\dof{p_i}} 
  \]
 where $2d+n$ is the number of edges of $\Gamma$ (see \cite{kontsevich}).
 
 The case that we will be most interested in is where $n=1$, in which case the above formulae simplify to give $\Omega=\omega$ with 
\be\label{eq:kontomega}
\omega=\sum_{1\leq i< j < m} \dof{l_i}\wedge\dof{l_j}
\ee
  where  $m$ is twice the number of  edges of $\Gamma$ and
  \be  \label{eq:volratio}
 \abs{\Omega^d}=\frac{4^dd!}{2^{g}}  \prod_{e}\abs{\dof{l(e)}}
  \ee
 where  the product  ranges over all but one edge of $\Gamma$.  Here we bring to the reader's attention that the form defined by \eqref{eq:kontomega} has a great deal of combinatorial symmetry (in particular, it is well-defined), some of which will implicitly be exploited in future calculations.

 The connection between the symplectic form $\omega$ and the tautological classes is given by the following
 \begin{thm}[Kontsevich \cite{kontsevich}]\label{thm:kont}
 The pullback of the form $\omega$ to $\Mb_{g,1}$ represents the class $\psi$: 
 \[
[ q^*\omega ] = \psi.
 \]
 where $q:\Mb_{g,1}\ra\Mcb_{g,1}$ is the composition of the quotient map $\Mb_{g,1}\ra\Mb'_{g,1}$ and the homeomorphism $\overline{\Psi}^{-1}:\Mb'_{g,1}\ra\Mcb_{g,1}$.
\end{thm}

The corresponding  theorem for $n>1$ is also valid  and provides a way of determining intersection products of the $\psi_i$ classes in $R^\bullet(\Mb_{g,n})$ by means of PL-integration over $\Mc_{g,n}$. 
In particular, for $n=1$ the theorem together with equation \eqref{eq:volratio} and the fact that the volume of a $1/2$-scaled $n$-simplex is $\frac{1}{2^n n!}$ leads directly to 
 \[
\int_{\Mb_{g,1}} \psi^{3g-2} =\int_{\Mc_{g,1}} \omega^{3g-2} =\frac{ (3g-2)!}{2^{g}(6g-4)! }\sum_\Gamma \frac{1}{\abs{\Aut(\Gamma)}}
\]
where the sum is over isomorphism classes of trivalent fatgraphs of type $(g,1)$ and the sign is deduced from the fact that all intersection products of $\psi$ classes are positive.  For $g=1$, one can immediately deduce that $\int_{\Mb_{1,1}} \psi =\frac{1}{24}$ as there is only one trivalent fatgraph of type $(1,1)$ and the automorphism group of this fatgraph has order 6.

\subsection{Combinatorial Cycles}
For $a\geq 0$, let $W_a\subset\Mgnc$ denote the locus of metric fatgraphs which have a  vertex of valence at least $2a+3$.   Witten was perhaps the first to suspect that the maximal cells in these loci can be endowed with a consistent orientation, thus they define locally-finite cycles representing codimension-$2a$  homology classes $[W_a]\in H^{\bullet-2a}_{lf}(\Mgnc)$. 

The loci $W_a$ are not quite orbifolds but rather are orbispaces with well defined fundamental classes.  The obstacle to being orbifolds is that the neighborhood in $W_a$ of a metric fatgraph with a vertex of  valence greater than $2a+3$ is not necessarily homeomorphic to the quotient of an open Euclidean ball.   For example, we have seen that near a $6$-valent vertex there are six ``directions'' which preserve a $5$-valent vertex.  Thus there are points in the cycle $W_1$ that have neighborhoods which look like three  sheets intersecting. 

Given a  perimeter value $p\in\oSim_{n-1}$, the slice $W_a(p)\subset\Mgnc(p)\homeo\Mgn$ can also be endowed with an orientation, thus $W_a(p)$ defines a locally finite homology class $[W_a(p)]\in H_{lf}^{\bullet-2a}(\Mgn)$.  Although the class $[W_a(p)]$ may depend on the value of $p$, for generic choices the class  will equal the class of the intersection $[W_a] \cdot [\Mgnc(p)]$ provided by the locally finite \Poin duality of $\Mgnc$.  Note that for $n=1$ there is only one (normalized) perimeter which is obviously generic.   By abuse of notation, when a generic slice is assumed, we shall denote $W_a(p)$ simply by $W_a$ and refer to it as the $a$th Witten cycle.
 
It was conjectured by Witten that the $W_a$  were \Poin duals to the  tautological $\kappa_a$ classes.
Witten's conjecture has since been proven and substantially generalized.  We refer to \cite{penner4}\cite{kontsevich}\cite{arbarello}\cite{igusa}\cite{igusa2}\cite{igusakleber}\cite{mondello} for the various perspectives on this conjecture and here only address that which is relevant to our discussion. In particular, while Witten's original conjecture concerned only the cycles $W_a$ on the open moduli, we shall be interested in the generalized result for $\Mgnb$. 

Kontsevich used  the symplectic form $\Omega$ to endow the $W_a$ with their orientations and consequentially showed that the Witten cycles could be naturally lifted to the partial compactification $\Mgncb$ where they could then be paired with products of the descendant $\psi_i$ classes.  However, since $\Mgncb$ is not an orbifold, the lack of \Poin duality makes it awkward to formalize what is meant by being dual to the $\kappa_a$ classes.  More recently, however,  Mondello has  shown by use of Looijenga's construction \cite{looijenga} that the Witten cycles in fact lift to cycles $\Wb_a$ in $\Mgnb$, where a full intersection theory exists.\footnote{Although this is in principle not needed for our result,  it does simplify the technical details considerably.  See \cite{mondello} for some taste of  the technical details that otherwise arise as well as a discussion  of the difficulty in lifting combinatorial cycles to $\Mgnb$.} 

We can now state the version of the theorem for $\Mgnb$.
\begin{thm}[Mondello \cite{mondello}] \label{thm:generalduality}
Under \Poin duality in $H^\bullet(\Mgnb)$ 
\[
\kappa_a\iso \frac{1}{2^{a+1}(2a+1)!!} \left([\overline{W}_a] + \delta_a \right)
\]
where $\delta_a$ is some boundary class. 
\end{thm}

In particular, we will make use of the concrete form of the theorem for $a=1$ which states
\begin{thm}[Mondello \cite{mondello}] \label{thm:duality}
Under \Poin duality, 
\[
\kappa_1\iso \frac{1}{12} ([\Wb_1]+ [\partial\Mgn]).
\]
\end{thm}

Immediately, we get as a corollary that
\be \label{eq:maincalcg1}
\int_{\Mb_{1,1}}\kappa_1= \frac{1}{12}\int_{\partial\M_{1,1}}1=
\frac{1}{24}
\ee
 as there no metric fatgraphs of type $(1,1)$ with a vertex of valence five, and $\partial\M_{1,1}$ is a single point with automorphism group of order two. 
 
 As a final remark, although the cycles $W_a$ represent algebraic classes, the actual loci are not themselves algebraic.  One can see this in the above  description of a neighborhood of a $6$-valent vertex as it gives an example of a real codimension one self-intersection of $W_1$.

\section{The Hyperelliptic Locus $\Hg^1$}
\label{sect:he}

In this section, we will describe the locus of hyperelliptic curves in terms of fatgraphs. We shall restrict to the case where $g\geq 2$, but almost everything will also hold for the elliptic case $g=1$.  There are several common definitions of hyperellipticity, and we shall be primarily interested in the one that states that a curve  $C$ is hyperelliptic if and only if $C$ has an automorphism  of order two  with $2g+2$ fixed points.  This automorphism is called the hyperelliptic involution and is denoted $\iota$.  A second related definition of hyperellipticity states that $C$ is hyperelliptic if and only if $C$ is a double covering of the Riemann sphere.  The involution $\iota$ then corresponds to the swapping of the two sheets of $\CP^1$, and the $2g+2$ fixed points are exactly the ramification points of this covering which can also be identified intrinsically as the \Weier points of the hyperelliptic curve $C$.

The collection of all hyperelliptic curves $\Hg$ forms an algebraic  subvariety of $\Mg$ of complex dimension $2g-1$.  As we are only interested in moduli spaces with $n>0$, instead of the loci of hyperelliptic curves, we shall rather consider the loci of hyperelliptic \Weier points $\Hg^1\subset\M_{g,1}$, which is a $2g+2$-fold ramified cover of $\Hg$ under $\pi_1$.   A point in $\Hg^1$ can be thought of as a pair $(C,x)$ with $C$ hyperelliptic and $x$ a \Weier point. The advantage of having $x$ so chosen is that the automorphism group $\Aut(C,x)$ of $(C,x)$ while not necessarily being  isomorphic to that of $C$ itself, is a subgroup of $\Aut(C)$ which contains the hyperelliptic involution $\iota$.   More generally, we shall also denote by $\He_{g,n}^m\subset\M_{g,m+n}$ the set of equivalence classes of hyperelliptic curves marked by $m$ \Weier points and $n$ arbitrary (distinct) points.  

\subsection{Hyperelliptic Fatgraphs}
Our goal now is to define  the locus of hyperelliptic metric fatgraphs $\Hgc\subset \Mc_{g,1}$ such that $\Psi(\Hgc)=\He_g^1$.  First we make the following preliminary definition:  a fixed point of a metric fatgraph automorphism $a\in\Aut(\Gamma,l)$ is a vertex, edge, or boundary cycle of $\Gamma$ which is taken to itself under $a$.  With this said, we make the 
\begin{defn}
A hyperelliptic fatgraph is a metric fatgraph $(\Gamma,l)$ which has an automorphism of order two with $2g+2$ fixed points.      We call this the hyperelliptic involution of $(\Gamma,l)$ and denote it also by $\iota$.  We also define the combinatorial hyperelliptic locus $\Hgc$ to be the locus of all hyperelliptic fatgraphs of type $(g,1)$ in $\Mc_{g,1}$.
\end{defn}

Given a metric fatgraph $(\Gamma,l)$ of type $(g,1)$, let $(C,x_1)=\Psi(\Gamma,l)$ be the corresponding  curve with marked point in $\M_{g,1}$.  Note that since $\Psi$ is a homeomorphism of orbifolds, $\Aut(\Gamma,l)$ is necessarily isomorphic to $\Aut(C,x_1)$.  What is more, one can easily prove that there is a bijection between the fixed points of an automorphism $a\in \Aut(\Gamma,l)$ and the fixed points of the corresponding automorphism $\Psi(a)\in\Aut(C,x_1)$.  
 In particular, a metric fatgraph of type $(g,1)$ will be hyperelliptic if and only if the corresponding marked curve is also.  In other words,  our definition is the appropriate one so that $\Psi(\Hgc)=\Hg^1$.

\begin{figure}[!h]
\psfrag{Hg}{$\Gamma_{H_2}$}
\psfrag{Hgp}{$\Gamma_{H_2}'$}
\begin{center}
\includegraphics[width=3.8in]{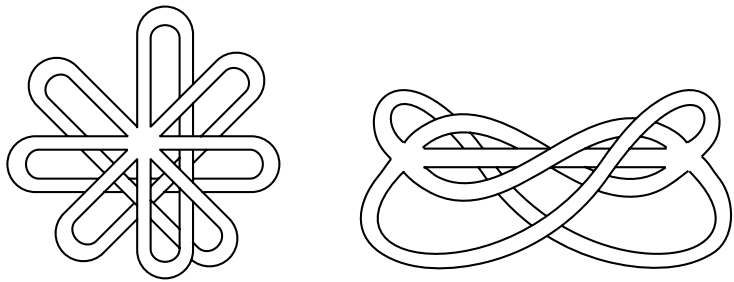}
\caption{Hyperelliptic ribbon graphs corresponding to minimal cells of $\Hgc$.}
\label{img:minhyperelliptic}
\end{center}
\end{figure}

\subsection{Minimal Cells of $\Hgc$}
Before giving a full description  of the locus $\Hgc$, we begin by identifying its lowest dimensional cells which are somewhat special and were first discussed by Penner in \cite{penner4}.

 Consider first  the  fatgraph  $\GHg$ which has only one vertex and $2g$ edges glued so that stubs of the same edge are opposite of each other in the cyclic ordering at the vertex.  Such a fatgraph has an automorphism group of order $4g$.  The first graph in figure \ref{img:minhyperelliptic} illustrates the case when $g=2$.  Notice that regardless of the metric put on this graph, there will always be an automorphism  of order two (rotation by $180^\circ$) which fixes all $2g$ edges, the one vertex, and the one boundary cycle.  Thus the corresponding orbi-simplex $\oSim_{\Gamma_{H_g}}/\Aut(\Gamma_{H_g})$ is contained in $\Hgc$.

Similarly, the second graph in figure \ref{img:minhyperelliptic} shows the $g=2$ case of a second type of fatgraph $\GHg'$ which has two vertices, $2g+1$ edges, and an automorphism group of order $2(2g+1)$.  Note that any edge collapse of  $\GHg'$ results in the graph $\GHg$, and that with any metric, there  is again an automorphism of order two with $2g+2$ fixed points.  Thus again, the orbi-simplex $\oSim_{\GHg'}/\Aut(\GHg')$ is contained in $\Hgc$.   

It can be proved that for no other (isomorphism class of) fatgraph $\Gamma$ is the full orbi-simplex $\oOrb$ contained in $\Hgc$; however, a quick dimension count shows that we have not described all of $\Hgc$.  In fact, we have only described the lowest dimensional cells of $\Hgc$, and we now turn towards  describing the remaining cells.  

\subsection{Trees and Spheres}
\label{sect:trees}
As hyperelliptic curves are all double covers of the Riemann sphere, they can be parametrized by their $2g+2$ branch points in $\CP^1$.  
Thus, it will be convenient to first give a  description of  the space $\M_{0,n}$ in terms of fatgraphs.  

Consider the slice of the genus zero combinatorial space $\Md_{0,n}(p)$ determined by the perimeter vector $p=(1,0,\ldots,0)$.  The cells of $\Md_{0,n}(p)$ are enumerated by fatgraphs of type $(0,1)$ with $n-1$ $\Delta$-labeled vertices.  As any fatgraph of type $(0,1)$ is necessarily a planar tree, we see that the maximal cells of  $\Md_{0,n}(p)$ are enumerated by trivalent trees with $n-1$ $\Delta$-labeled leaves (uni-valent vertices).  It is easy to check that the number of edges of such a tree is $2n-5$, and there is one constraint on the possible metrics, in agreement with the complex dimension of $\M_{0,n}$ being $n-3$.

It is straightforward to show that the relation $[\omega_1]=\Psi^*(\psi_1)$ extends to the partially compactified space  $\Md_{0,n}(p)$ for this particular $p$ (all that is needed is for $p_1\neq0$), and that  Kontsevich's symplectic form for this slice is just $\Omega=\omega_1$.  The volume determined by $\Omega$ on a simplex is provided by the following 

 \begin{lemma} \label{volume}
 Let $T$ be a tree with all vertices of odd valence and $2m+1$ distinct edges $e_i$.  Then on  the simplex $\oSim_T$ corresponding to the slice  $\rho_0=(1,0,\dotsc,0)$, we have 
 \[
 \abs{\omega_1^m}=4^mm! \prod_{i=1}^{2m}\abs{\dof{l(e_i)}}.
 \]
 \end{lemma}
 \begin{proof}
 This is an immediate corollary of Kontsevich's result \eqref{eq:volratio}; however, it is also easy to prove the result by induction on $m$ directly from the definition \eqref{eq:kontomega} of the form $\omega$.  The base case $m=0$ corresponds to a tree with a single edge and is obvious, while 
 for $m >0$, we note that any tree  of odd valence with three or more leaves must have a pair of consecutive leaves which meet at a common vertex.  By removing these two leaves, we obtain a new tree $T'$ which satisfies the induction hypothesis for $m-1$. We leave the remaining details to the reader.
 \end{proof}
 
 As a corollary, we can verify the integral
  \[
 \int_{\Mb_{0,n}}\psi^{n-3}=\sum_T \int_{\oSim_T}\omega^{n-3}=(n-2)! C_{n-3}\frac{(n-3)!}{(2n-6)!}=1
 \]
 where the sum is over trivalent trees with $n-1$ labeled leaves, and $C_{n-3}$ is the Catalan number which counts the number of rooted trivalent trees with $n-3$ trivalent (i.e.\ non-leaf) vertices.\footnote{This result actually requires a common extension of the two compactifications $\Mgnd$ and $\Mgncb$, and we refer to \cite{looijenga} or \cite{mondello} for more details.}  This appearance of the Catalan numbers is a simple illustration of  a more general phenomenon which we will later encounter in the main calculation of the paper.  Finally, we note that it is instructive to apply similar techniques and  theorem \ref{thm:duality} to compute the integral $\int_{\Mb_{0,n}} \psi^{n-4}\kappa_1$, and we recommend the interested reader to do so.  

\subsection{Cell Decomposition of $\Hgc$}
\begin{figure}[!h]
\begin{center}
\includegraphics[width=1.4in]{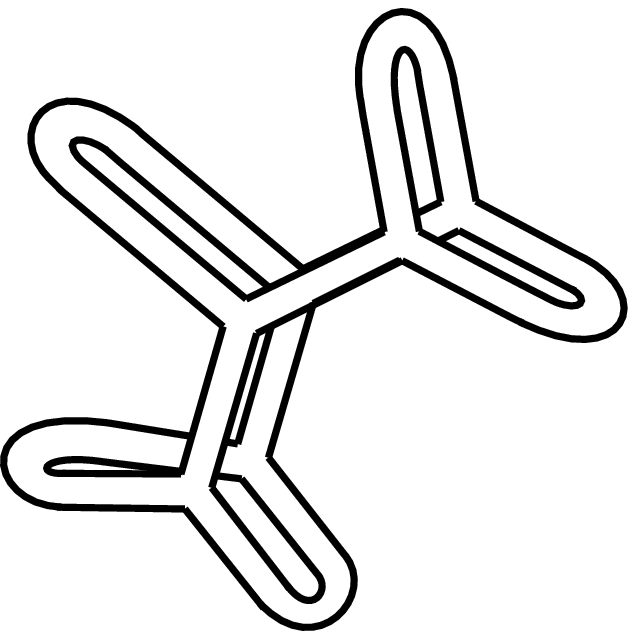} \hspace{.6 in}
\includegraphics[width=1.9in]{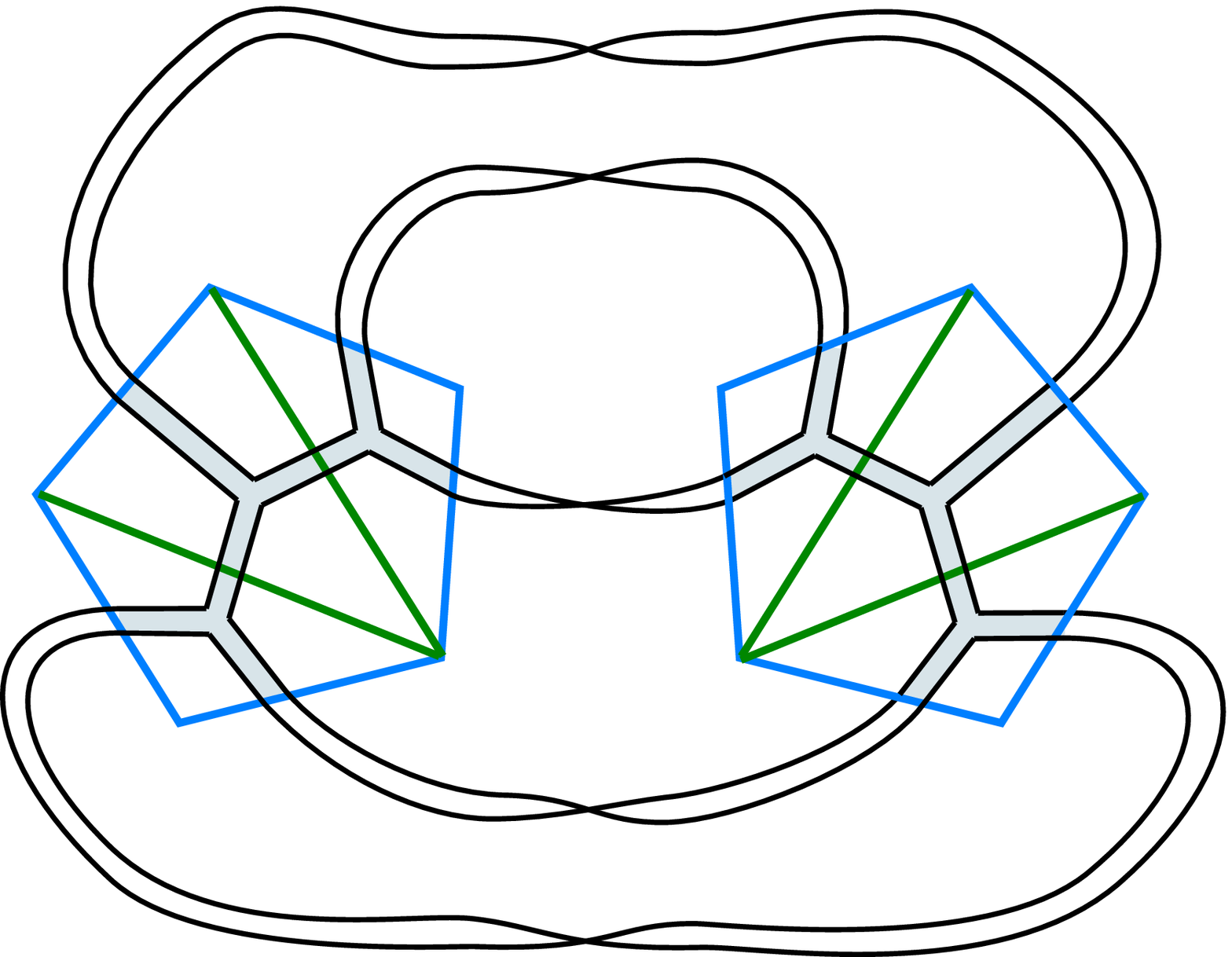}
\caption{Equivalent trivalent hyperelliptic  ribbon graphs.}
\label{img:hyperelliptic}
\end{center}
\end{figure}

Consider now a trivalent hyperelliptic fatgraph $(\Gamma,l)$ of type $(g,1)$ (the existence of which will soon be obvious).   Two topologically equivalent ribbon graph depictions of a genus 2 trivalent hyperelliptic fatgraph are given in figure \ref{img:hyperelliptic}.  As $\Gamma$ is trivalent, we immediately see that no vertex can be fixed by the hyperelliptic involution, thus there are $2g+1$ fixed edges.  Now, we construct a new disconnected metric fatgraph from $(\Gamma,l)$ by ``cutting'' each fixed edge in two equal parts, thus resulting in $4g+2$ new univalent vertices.  It is easy to prove by induction on $g$ that the resulting graph consists of two identical components, each of which is a trivalent tree with $2g+1$ leaves. 
 
\begin{figure}[!h]
\begin{center}
\psfrag{ep1}{$\varep_1$}
\psfrag{ep2}{$\varep_2$}
\psfrag{ep3}{$\varep_3$}
\psfrag{ep4}{$\varep_4$}
\psfrag{ep1p}{$\varep'_1$}
\psfrag{ep2p}{$\varep'_2$}
\psfrag{ep3p}{$\varep'_3$}
\psfrag{ep4p}{$\varep'_4$}
\psfrag{e1}{$e_1$}
\psfrag{e2}{$e_2$}
\psfrag{e3}{$e_3$}
\psfrag{e4}{$e_4$}
\psfrag{e5}{$e_5$}
\psfrag{e6}{$e_6$}
\psfrag{e7}{$e_7$}
\includegraphics[width=3.3in]{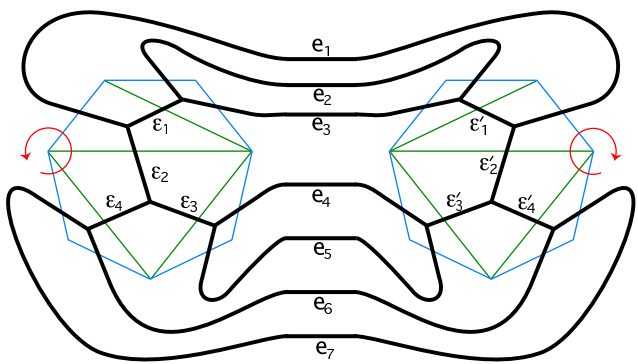}
\caption{Edges of $\Gamma_T$ labeled.}
\label{img:hyperlabelled}
\end{center}
\end{figure}

\begin{figure}[!h]
\begin{center}
\psfrag{Dh}[][]{$\oSim_{H_{T}}$}
\psfrag{Dg}[][]{$\Delta_{\Gamma_T}$}
\includegraphics[width=1.6in]{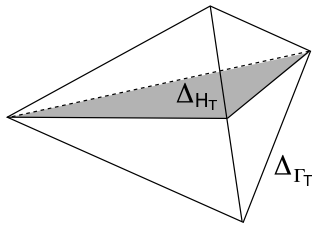} 
\caption{Hyperelliptic cell $\oSim_{H_T}$ in first barycentric subdivision of $\Delta_{\Gamma_T}$.}
\label{img:hecell}
\end{center}
\end{figure}

Conversely, given a trivalent (normalized) metric tree $(T,l^T)$ with $2g+1$ leaves, we can make an identical copy,  glue  corresponding leaves, and re-scale to obtain a hyperelliptic fatgraph $(\Gamma_T,l)\in\Hgc$.  To explicitly define the metric $l$ on $\Gamma_T$ in terms of  $l^T$, we first label the edges of the tree $T$ by $\{e^T_1,\ldots,e^T_{2g+1},\varep^T_1,\varep^T_{2g-1}\}$ where the edges $e^T_i$ correspond to the leaves and the $\varep^T_j$ correspond to the remaining edges.   Next we label the corresponding edges of $\Gamma_T$ by $\{e_1,\ldots,e_{2g+1},\varep_1,\ldots,\varep_{2g-1},\varep'_1,\ldots,\varep'_{2g-1}\}$ as shown in figure \ref{img:hyperlabelled}.  Finally we define the metric $l$ on $\Gamma_T$ by 
 \be
 \label{eq:treetohyper}
 l(e_i)= l^T(e^T_i), \hspace{.5 in}
 l(\varep_j)=l(\varep'_j)= \frac{1}{2}\cdot l^T(\varep^T_j).
 \ee
  Equations \eqref{eq:treetohyper} define an inclusion map  $h:\oSim_T \hookrightarrow \oSim_{\Gamma_T}$, and we denote the image by $\oSim_{H_T}$, as depicted in figure \ref{img:hecell}.   
 
We now observe that the above map  extends in a well defined way  to the faces of $\oSim_T$ contained in $\Md_{0,2g+2}(p)$.  In particular, given any $\Delta$-labeled tree $T'$ of type $(0,2g+2)$, we can double it to obtain a genus $g$ hyperelliptic fatgraph  where any $\Delta$-labeled vertices of valence $v_i>1$ are glued to their doubles in any of the $v_i$ possible symmetric ways.  Conversely, given an arbitrary hyperelliptic fatgraph, we can obtain two identical trees by splitting symmetrically along the $2g+1$ non-boundary fixed points, whether they are edges or vertices (note that this is well defined).  
 
 The automorphism group $\iota$ acts trivially on $\oSim_{H_{T'}}$  while $\Aut(\Gamma_{T'})/\iota$  will always act faithfully and in a way compatible with the gluing of cells.  Thus we get the desired cell decomposition of $\Hgc$, and we summarize with the
 \begin{lemma} The locus of hyperelliptic \Weier points has a combinatorial cell decomposition given by 
 \be\label{eq:hecelldecomp}
 \Psi^{-1}(\Hg^1) =  \Hgc=\coprod_T \oSim_{H_T}/\Aut(\Gamma_T)
  \ee
 where the union is taken over equivalence classes of $\Delta$-labeled trees $T$ of type $(0,2g+2)$, and two trees are considered equivalent if they are isomorphic as fatgraphs and their sets of (unordered) $\Delta$-labeled vertices coincide.
 \end{lemma}
   In fact, we have implicitly just described the PL-orbifold covering map,
 \[
 \Md_{0,2g+2}(p) \ra \Hgc
 \]
 for $p=(1,0,\ldots,0)$.  

For one familiar with Strebel differentials, the above correspondence between hyperelliptic graphs and planar trees will be more or less obvious in the conformal setting as a Strebel differential and its corresponding critical trajectories on $\CP^1$ will both be pulled back to the hyperelliptic curve under the double covering.  

The following is the main technical lemma for the computations which will follow in the next section.
\begin{lemma}  \label{hgvol}
Let $(\Gamma_T,l)$ be a hyperelliptic graph  corresponding to a tree $T$ which has only odd valence vertices and $2d+1$ edges.  Then,
\[
\int_{\oSim_{H_T}}\omega^d = \frac{d!}{2^d(2d)!}.
\]
\end{lemma}
\begin{proof}
This lemma essentially follows from lemma \ref{volume} and the fact that Kontsevich's two-form on the cells $\oSim_{H_T}$ and $\oSim_T$ is related by
\be\label{eq:omegafact}
h^*\omega_{H_T}=\frac{1}{2}\cdot \omega_T.
\ee
This fact follows from two observations.  To see these, it is useful to first ``split''  each fixed edge of  $\Gamma_T$ by inserting a bivalent vertex.   Next,  identify the edges of each ``half of'' $\Gamma_T$ with the corresponding edges of $T$.   From this, one can see that the boundary cycle of $\Gamma_T$ is essentially equal to that of $T$ repeated twice.  Thus by a careful application of the definition \eqref{eq:kontomega}, the form $\omega_{H_T}$ looks identical to $2\cdot\omega_T$, except that it is expressed in terms of the length coordinates for $\oSim_{H_T}$.  As these length coordinates  are all half scaled with respect to those of $\oSim_T$ (C.f. \eqref{eq:treetohyper}), the $2$-form must be scaled by a factor of $\frac{1}{4}$.  Thus the two factors combine to give \eqref{eq:omegafact}.

Finally, lemma \ref{volume} immediately yields 
\[
\int_{\oSim_T}\omega_T^d=4^d d!\int_{\oSim_T} \prod_{i=1}^{2d}\abs{\dof{l}(e_i)}=\frac{d!}{(2d)!},
\]
which together with equation \eqref{eq:omegafact} gives the result, where again the sign has been determined by the positivity of intersection numbers.
\end{proof}

In particular, the above lemma  shows that the locus $\Hgc$ is symplectic with respect to $\Omega=\omega$ as the volume form $\omega^d$ is non-degenerate. 

\section{Calculations}
\label{sect:calc}

As the locus of hyperelliptic points $\Hg^1$ is an algebraic subvariety of $\M_{g,1}$, it has a natural closure $\Hgb^1$ in $\Mb_{g,1}$.   Thus it makes sense to integrate tautological classes over this locus.

As a first example, we can compute the integral of $\psi^{2g-1}$ over the locus $\Hgb^1$ by employing the integration scheme  \eqref{eq:intscheme} with respect to the cell decomposition \eqref{eq:hecelldecomp} and the representation $\psi=[\omega]$,
\[
\int_{\Hgb^1}\psi^{2g-1} =\int_{\Hgc} \omega^{2g-1}=\sum_T \frac{1}{\abs{\Aut(\Gamma_T)}} \int_{\oSim_{\Gamma_T}}\omega^{2g-1}.
\]
Using lemma \ref{hgvol} and the Catalan numbers \eqref{eq:catalan} to count (in the orbifold sense) the maximal hyperelliptic cells corresponding to  trivalent trees with $2g+1$ leaves, we get (C.f.\ \cite{faberp})
\be\label{eq:hevol}
\int_{\Hgb}\psi^{2g-1} =\frac{C_{2g-1} }{2(2g+1)}\cdot  \frac{(2g-1)!}{2^{2g-1}(4g-2)!}=\frac{1}{2^{2g}(2g+1)!}.
\ee

The main theorem of this paper is the calculation of a similar yet more difficult integral.
\begin{mainthm}\label{thm:maincalc}
For $g\geq1$, the integral of tautological classes $\kappa_1\psi^{2g-2}$ over the locus of hyperelliptic \Weier points is given by 
\[
\int_{\Hgb^1}\kappa_1\psi^{2g-2}=\frac{(2g-1)^2}{2^{2g}(2g+1)!}.
\] 
\end{mainthm}

Note that for the elliptic case $g=1$, the result  has already been shown by equation \eqref{eq:maincalcg1}, thus  we shall assume $g>1$. The proof now proceeds in several steps.

First, we use theorem \ref{thm:duality} and the \Poin duality of $\Mb_{g,1}$  to rewrite the integral as
\be \label{eq:maincalc}
\int_{\Hgb^1}\kappa_1\psi^{2g-2}   =   \frac{1}{12} \left[   \int_{[\Wb_1]\cdot [\Hgb^1]}\psi^{2g-2} + \int_{[\partial\Mb_{g,1}]\cdot[\Hgb^1]}\psi^{2g-2}     \right].
\ee
The heart of the calculation is the evaluation of the first integral on the right hand side of \eqref{eq:maincalc}.   Assuming for the moment that  the loci $\Wb_1$ and $\Hgb^1$ intersect transversely,
this term is calculated by integrating the form $\omega^{2g-2}$ representing $\psi^{2g-2}$ over maximal cells of each component of the intersection  $W_1\cap\Hgc$ in  $\Mc_{g,1}$  and adding these contributions with appropriate multiplicities.  The integral over each cell will be seen to give the same volume, thus the total integral will be reduced to a fatgraph counting argument much as was done in \eqref{eq:hevol}.

In sections \ref{sect:intw1hg} and \ref{sect:count} we will compute the first integral of equation \eqref{eq:maincalc} while  the second integral is addressed in section \ref{sect:heboundary}.  In section \ref{sect:corr} we use the above result to derive a corollary about Hodge integrals. 

\subsection{The Intersection $W_1\cap\Hgc$}
\label{sect:intw1hg}
By definition,  the  locus  $W_1\cap\Hgc$ consists of hyperelliptic fatgraphs which have a vertex of valence at least 5.  One can check that the maximal cells in this locus have (real) dimension $4g-4$ and come in two types, and the union of all (closures of) cells of each type forms a component of the intersection $W_1\cap\Hgc$.  Moreover, an application of lemma \ref{hgvol} will show that each component is symplectic with respect to $\Omega$, thus inherits a natural orientation compatible with those of $W_1$ and $\Hgc$. 
 We shall investigate each component in turn.  

\begin{figure}[!h]
\begin{center} 
\includegraphics[width=2.0in]{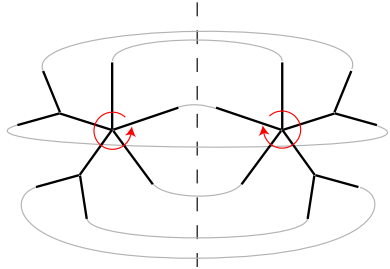}
\caption{Hyperelliptic fatgraph in first component of $W_1\cap\Hgc$.}
\label{img:twotrees5val}
\end{center}
\end{figure}

The first type of cell in $W_1\cap\Hgc$ is enumerated by trees with $2g+1$ leaves which are trivalent except for a single vertex of valence five.  A   hyperelliptic fatgraph in one of these  cells will thus necessarily have two vertices of valence five, as the example depicted in   figure \ref{img:twotrees5val} shows for the case $g=3$.   As a result of there being not one but two $5$-valent vertices, the multiplicity of this locus is two.  We note that this factor of two in some ways cancels the hyperelliptic symmetry factor in the orbifold fundamental class of this component.\footnote{More accurately, each component of $W_1$ corresponds to a choice of $5$-valent vertex.  However, there is essentially only one choice, as both are isomorphic under $\iota$.  This ``loss'' of a factor of two is balanced by the fact that once such a choice is made, the (decorated)  graph is no longer invariant under the hyperelliptic involution, thus it also ``loses'' a symmetry factor of one half.}
  
\begin{figure}[!h]
\begin{center}
\includegraphics[width=3.0in]{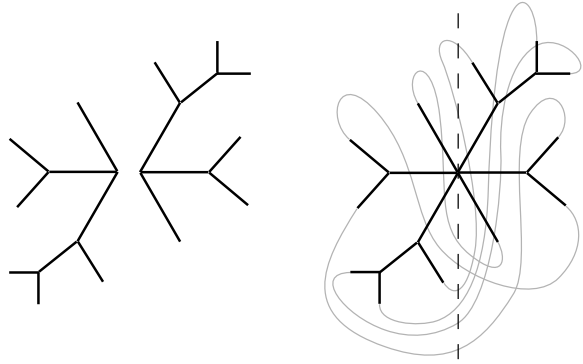}
\caption{Hyperelliptic fatgraph in second component of $W_1\cap\Hgc$.}
\label{img:twotrees6valg4}
\end{center}
\end{figure}

 The second type of cell is enumerated by trivalent trees with $2g$ leaves and one $\Delta$-labeled trivalent vertex.  A hyperelliptic fatgraph $\Gamma$ in one of  these cells then will have a single fixed vertex of valence six, and an example for $g=4$  is depicted in   figure \ref{img:twotrees6valg4}.   We have already mentioned that at a $6$-valent vertex, three `sheets' of $W_1$ intersect,  thus we suspect that the multiplicity of this locus should be 3. In order to concretely determine this, we need to show  that   each of the three sheets of $W_1$ intersects $\Hgc$ transversely.  We note here that some difficulty in determining this multiplicity is caused by the orbifold nature of the cycles $W_1$ and $\Hg^1$.  In particular, the symmetry of the hyperelliptic involution acts on these cycles, and to determine the proper multiplicity of the intersection, we will analyze it in the lift to a Euclidean neighborhood  in \Teich space.   
 
\begin{figure}[!h]
\begin{center}
\includegraphics[width=4.0in]{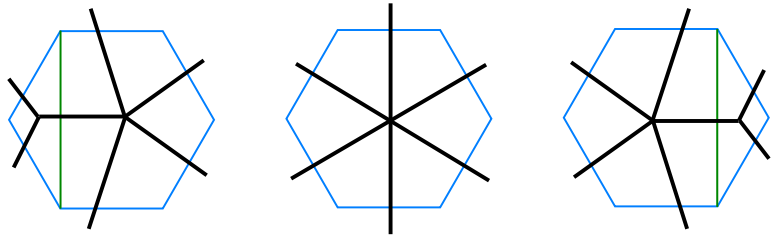}
\caption{Two opposite Witten expansions of a $6$-valent vertex.}
\label{img:oppositeorientations}
\end{center}
\end{figure}
 
 Recall that there are nine ways to expand a $6$-valent vertex to produce one new edge, corresponding to the nine ways of drawing a ray between vertices of a regular hexagon.   Six of these expansions  correspond to  cells lying in $W_1$ and can be put in 1-to-1 correspondence with the vertices of the hexagon (corresponding to which vertex is ``cut off').   Two \emph{opposite} expansions are depicted in figure \ref{img:oppositeorientations}.   As $W_1$ forms a cycle,  the orientations  on the six  cells  must produce canceling contributions at their common face.  In fact, a  more careful analysis  (perhaps easiest in terms of Conant and Vogtmann's definition of orientation \cite{vogtmann}) shows that any  two \emph{adjacent} or \emph{opposite} expansions (corresponding to adjacent or opposite vertices) must contribute canceling orientations.  Thus it is natural to construct the cycle $W_1$ by pairing  cells of opposite expansions to form the sheets of $W_1$.
 
\begin{figure}[!h]
\begin{center}
\includegraphics[width=4.0in]{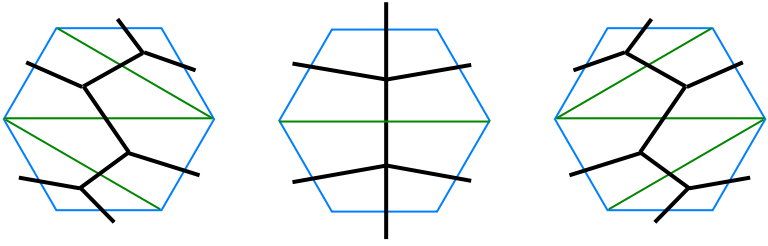}
\caption{Hyperelliptic expansions of a $6$-valent vertex.}
\label{img:heexpansions}
\end{center}
\end{figure}

 The three remaining one-edge expansions which do not correspond to cells in $W_1$ result in a graph which has two $4$-valent vertices and a $180^\circ$ rotational symmetry.    As a result, the corresponding cells retain the symmetry of the hyperelliptic involution, thus lie in $\Hgc$.  However, there are other expansions which introduce two or three edges which (with appropriately symmetric metric assignments) also correspond to cells  lying in $\Hgc$.  See figure \ref{img:heexpansions}.  
 
 One can check that (points in the interiors of) the paired opposite cells of $W_1$ cannot be connected (in the Euclidean neighborhood) without passing through a fatgraph with the hyperelliptic symmetry.  Moreover, a quick check shows that at the $6$-valent vertex, the Witten cycle $W_1$ meets $\Hgc$ in complimentary dimension.  
 Thus we see that the intersection $W_1\cap\Hgc$ is in fact transverse and has multiplicity  three as predicted, a contribution of one being provided by each sheet of paired opposite cells.  Note that the signs of intersections are automatically determined from symplectic structures of the cycles involved  (see below).  
We summarize the above discussion in figure \ref{img:HgW1intersection}, where we have  depicted the  intersection  together with the associahedron $K_5$.   

\begin{figure}[!h]
\psfrag{W}{$W_1$}
\psfrag{Hg}{$\Hgc$}
\begin{center}
\includegraphics[width=3.0in]{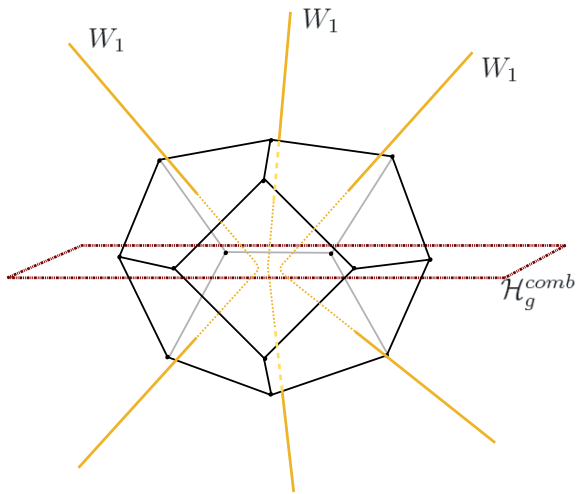}
\caption{The intersection of lifts of $W_1$ and  $\Hgc$ in the associahedron $K_5$.}
\label{img:HgW1intersection}
\end{center}
\end{figure}

\subsection{Enumeration of Cells in $W_1\cap\Hgc$}
\label{sect:count}
Now that we have determined the multiplicities of the two components of $W_1\cap\Hgc$, we can rewrite the first integral on the RHS of \eqref{eq:maincalc}  as  a weighted sum over graphs 
\be \label{eq:w1hgcsum}
\int_{[\Wb_1]\cdot[\Hgb^1]}\psi^{2g-2} = 2\sum_{T_1}\frac{1}{\abs{\Aut(\Gamma_{T_1})}}\int_{\oSim_{H_{T_1}}}\!\omega^{2g-2} + 3\sum_{T_2}\frac{1}{\abs{\Aut(\Gamma_{T_2})}}\int_{\oSim_{H_{T_2}}}\!\omega^{2g-2} 
\ee
where the first sum is over trees $T_1$ of the first type with a $5$-valent vertex  and the second sum is over trees $T_2$ of the second type with a $\Delta$-labeled trivalent vertex.  

As both types of trees $T_1$ and $T_2$ have $4g-3$ edges, by lemma \ref{hgvol} the symplectic volume of the corresponding cells $\oSim_{T_1}$ and $\oSim_{T_2}$ will both be 
\[
\int_{\oSim_{H_T}}\omega^{2g-2}= \frac{(2g-2)!}{2^{2g-2}(4g-4)!}.
\]
Consequentially,  this also verifies the non-degeneracy of the symplectic form $\Omega$ on $W_1\cap\Hgc$, and thus inherently accounts for the proper orientation of the components of the intersection.

As a result, the integral \eqref{eq:w1hgcsum} is reduced to two orbifold counting arguments,
\[
\int_{[\Wb_1]\cdot[\Hgb^1]}\psi^{2g-2} = \frac{(2g-2)!}{2^{2g-2}(4g-4)!}\left[ 2\sum_{T_1}\frac{1}{\abs{\Aut(\Gamma_{T_1})}}+  3\sum_{T_2}\frac{1}{\abs{\Aut(\Gamma_{T_2})}}\right] 
\]
 which we now address. First, we count the number of cells arising in the first component.  These cells are enumerated by trees with one vertex of valence 5 and $2g-4$ trivalent ones.   As we have seen, the number of rooted trees of this type is given by the generalized Catalan number $C_{5,2g-3}$.  As we rather wish to count the number of unrooted trees in the orbifold sense, we must divide by the number of leaves $2g+1$ as well as a factor of $2$ for the hyperelliptic involution.  Thus get that 
\be \label{eq:count1}
\sum_{T_1}\frac{1}{\abs{\Aut(\Gamma_{T_1})}}=\frac{C_{5,2g-3}}{2(2g+1)}.
\ee

For the second component, the cells are enumerated by trivalent trees with $2g$ leaves and one of its $2g-2$ trivalent vertices specially marked.  The number of such rooted trees  is just the product of the ordinary Catalan number $C_{2g-2}$ and the number of trivalent vertices $2g-2$.  To get the desired orbifold sum we must also divide by the number of leaves $2g$ and the order of the hyperelliptic involution $2$, thus getting
\be \label{eq:count2}
\sum_{T_2}\frac{1}{\abs{\Aut(\Gamma_{T_2})}}=\frac{(g-1)C_{2g-2}}{2g}.
\ee

Using equations \eqref{eq:count1}, \eqref{eq:count2}, \eqref{eq:catalan},  and \eqref{eq:gencatalan}, the integral \eqref{eq:w1hgcsum} can be evaluated to be 
\begin{multline} \label{eq:w1hgvol}
\int_{W_1\cap\Hgc} \omega^{2g-2}=  \frac{1}{2^{2g-2}(2g+1)!}\Big(  (2g-2)(2g-3) + 3(2g+1)(g-1)  \Big)\\
=\frac{10g^2-13g+3}{2^{2g-2}(2g+1)!}.
\end{multline}

\subsection{The Intersection $\Hgb^1\cap\partial\M_{g,1}$}
\label{sect:heboundary}
To complete our calculation, we need to determine the value of the integral
\[
\int_{[\partial\Mb_{g,1}]\cdot[\Hgb^1]}\psi^{2g-2}=\int_{\partial\Hg^1}\psi^{2g-2}
\]
where $\partial\Hg^1$ is the locus of the intersection $\Hgb^1\cap \partial \M_{g,1}$.

Each stable hyperelliptic \Weier point in $\partial\Hg^1$ lies on a curve with at least one node (although it is the limit of    \Weier points on smooth hyperelliptic curves).  Thus  $\partial\Hg^1$ has several irreducible components corresponding to how the (generically) one node breaks the curve into hyperelliptic components of different genera.  
 More precisely,
$\partial\Hg^1$ decomposes as the union of images under the  hyperelliptic  gluing morphisms 
\be \label{eq:heboundary1}
\delta^\He_{h,g-h}:\Hb_h^1\times\Hb_{g-h}^2\ra\Hb_{g}^1
\ee
 for $h=1,\ldots,g-1$ and 
 \be \label{eq:heboundary2}
 \delta^\He_{g-1}:\Hb_{g-1,1}^1\ra\Hb_{g}^1
 \ee
 where the first map identifies marked hyperelliptic \Weier points of different curves to form a node, and the second map identifies the non-\Weier point of a hyperelliptic curve with its image under the hyperelliptic involution.
 
 As we are only interested in the value of the integral $\int_{\partial\Hg^1}\psi^{2g-2}$, we shall only consider components of $\partial\Hg^1$ where this integral does not vanish. 
  For this to happen, $\psi^{2g-2}$ must not be identically zero, thus the \Weier  point  must  lie on a (possibly nodal) hyperelliptic  component with  at least $2g-2$ (complex) moduli.   This occurs  only on the  image of the map \eqref{eq:heboundary2}, and we may identify this component of $\partial\Hg^1$ with ``half of"  the universal curve $\Hb_{g-1,1}^1$ over $\Hb_{g-1}^1$.  The integral can now be evaluated as 
 \be\label{eq:hedegen}
\int_{\partial\Hg^1}\psi^{2g-2}=\frac{1}{2}\int_{\Hb^1_{g-1,1}}\psi^{2g-2}=\frac{1}{2}\int_{\Hb^1_{g-1}}\psi^{2g-3}=\frac{1}{2^{2g-1}(2g-1)!},
\ee
where the second equality follows from the string equation.

\begin{figure}[!h]
\psfrag{one}[][]{$1$}
\begin{center}
\includegraphics[width=4.8in]{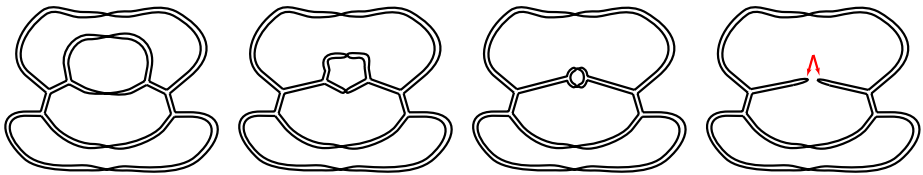}
\caption{Pinching of a genus $2$ hyperelliptic fatgraph, resulting in a stable hyperelliptic  fatgraph with genus $1$ node.}
\label{img:stablehe}
\end{center}
\end{figure}

We note that the above calculation can be performed alternatively by use of Kontsevich's stable fatgraphs.  In particular, in $\Mcb_{g,1}$, the images of components of $\partial\Hg^1$ have varying dimension, and the one of maximal dimension is precisely the one  corresponding to the map \eqref{eq:heboundary2}.  The stable graphs corresponding to the maximal cells of this component are 
easily visualized as degenerations of trivalent ones as in figure \ref{img:stablehe}.  As such, the cells are enumerated by rooted trivalent trees with $2g$ leaves, where the root corresponds to the node of the stable hyperelliptic curve.  One can then derive equation \eqref{eq:hedegen} much as was done for equation \eqref{eq:hevol}.
 
Putting the contributions \eqref{eq:w1hgvol} and \eqref{eq:hedegen} together, we can finally evaluate the integral \eqref{eq:maincalc} 
\begin{multline}\label{final}
\int_{\Hgb^1}\kappa_1\psi^{2g-2}  =  \frac{1}{12}\left[  \frac{1}{2^{2g-2}(2g+1)!}\Big( 10g^2-13g+3+ g(2g+1) \Big) \right]  \\
 =  \frac{(2g-1)^2}{2^{2g}(2g+1)!},
\end{multline}
 thus proving  theorem \ref{thm:maincalc}.

\subsection{Hodge Integrals}
\label{sect:corr} 
We now use the result of the last section to derive as a corollary a relation between Hodge integrals in the tautological ring of $\Mgnb$.  In particular we shall prove the following.
\begin{cor} \label{cor:cor}
For $g\geq2$, the following relation of Hodge integrals holds
\be
\sum_{i=0}^{g-1} \int_{\Mb_{g,1}} \left[ (-1)^{i}(2^{g-i}-1)\psi^{g-i-1}\lambda_{i} \right]\kappa_1\psi^{2g-2} =  \frac{14g^2 -11g +3  }{3\cdot{}2^{2g}(2g+1)!}.
\ee
\end{cor}

The derivation of this corollary relies on the  correspondence between the locus of hyperelliptic \Weier points and tautological classes, which on the open moduli is well known and  goes back to  Mumford \cite{mumford}.  It states that on $\M_{g,1}$, 
\be \label{eq:porteushe}
[\Hg^1]\underset{g-1}{=}   \frac{ 1 - \lambda + \lambda^2 - \dotsm  \pm \lambda^{g-1}}{ (1-\psi)(1-2\psi)},
\ee
where the equality holds only for the degree  $g-1$  terms.
Mumford's original derivation  relied on an application of the Grothendieck-Riemann-Roch formula; however, the relation can also be obtained by means of Porteus' formula, which  has the advantage that it can be partially extended to the compactified moduli space $\Mb_{g,1}$.

Porteus' formula (see \cite{harris}) gives a relation between the class of the degeneracy locus $X$ of a bundle map $E\ra F$ and the characteristic classes of  $E$ and $F$, under the condition that the locus $X$ is of the expected dimension.  In the case of rank one  degeneracy,  Porteus' formula takes the simple form
\be \label{eq:porteus}
[X]\underset{d}{=} \frac{c(E^*)}{c(F^*)}
\ee
where $c$ denotes the total Chern class, $E^*$ is the dual bundle of $E$, and $d$ is the codimension of $X$.

For the case at hand, the bundle map that is relevant  is the fiberwise  evaluation map 
\be  \label{degenerate}
H^0(C,\omega_C)\ra H^0(C,\omega_C/\omega_C(-2x))
\ee
for the Hodge bundle $\mathbb{E}$ and the expected dimension of degeneracy is $2g-1$ (codimension $g-1$).  By the classical Riemann-Roch theorem, the degeneracy of this map at $(C,x)$ for $C$ a smooth curve is equivalent to the existence of a meromorphic function on $C$ holomorphic except for a double pole at $x$.   In other words, the degeneracy locus in $\M_{g,1}$ exactly the locus of  hyperelliptic \Weier points and Porteus' formula \eqref{eq:porteus} gives \eqref{eq:porteushe}.

 For singular curves $(C,x)$, however, the  the above evaluation map will be degenerate whenever the point $x$ lies on a hyperelliptic, elliptic, or rational component of $C$, with $C$ itself not necessarily (the limit of) hyperelliptic.  One can check that this locus has components of varying dimension (both greater and lesser than expected), thus a direct application of Porteus' formula is not valid.  However, a more careful statement of Porteus' formula states that (see \cite{fulton}) 
 \[
[X]\underset{g-1}{=}    \frac{ 1 - \lambda + \lambda^2 - \dotsm  \pm \lambda^{g-1}}{ (1-\psi)(1-2\psi)}.
\]
where $X$ is some subvariety of dimension $2g-1$ contained in the degeneracy locus such that $X$ contains all irreducible components of the correct dimension.  In particular, $X=\Hgb^1 + Z$ with $Z\subset\partial\Mgnb$.

We now turn towards determining $Z$ and shall only be interested in components $Y\subset Z$ such that the integral $\int_Y \kappa_1\psi^{2g-2}$ does not vanish  (c.f.\ \cite{faberp}).   Note  that for this to happen, the single marked point  must  lie on a hyperelliptic (or elliptic or rational) component with  at least $2g-2$ (complex) moduli.  Also note that the presence of the $\kappa_1$ class makes this locus not describable solely in terms of the space $\Mcb_{g,1}$, unlike the case for the integral  \eqref {eq:hedegen}.

Recall that the boundary $\partial\M_{g,1}$ decomposes as the union of images under the  gluing morphisms 
\be \label{eq:boundary1}
\delta^1_{h,g-h}:\Mb_{h,1}\times\Mb_{g-h,2}\ra\Mb_{g,1}
\ee
 for $h=1,\ldots,g-1$ and 
 \be \label{eq:boundary2}
 \delta^1_{g-1}:\Mb_{g-1,3}\ra\Mb_{g,1}.
 \ee
 
 The components of $Z$ in the image of $\delta^1_{g-1}$ can be disregarded as the only one on which $\psi^{2g-2}$ does not vanish is  the locus $\delta^1_{g-1}(\Mb_{g-1,3})\cap\Hb_{g}^1$, which  we have seen   has dimension $2g-2$, one less than the required $2g-1$.

For components  in the images of the $\delta^1_{h,g-h}$, the only possibility of non-vanishing $\psi^{2g-2}$ occurs when $h=1$ and $g>1$ on the locus  $Y=\delta^1_{1,g-1}(\Mb_{1,1}\times\Hb_{g-1,1}^1)$.  Although curves in $Y$ may not necessarily be hyperelliptic, they can be described as hyperelliptic components with elliptic tails attached at a non-\Weier node, and $Y$ is easily seen to have the correct dimension $2g-1$.  

Thus by Porteus' formula, we are able to arrive at the following for $g\geq 2$
\[
\int_{\Mb_{g,1}}  \frac{ 1 - \lambda + \lambda^2 - \dotsm  \pm \lambda^{g-1}}{ (1-\psi)(1-2\psi)} \kappa_1\psi^{2g-2} = \int_{\Hgb^1}\kappa_1\psi^{2g-2} + \int_{Y}\kappa_1\psi^{2g-2}.
\]

The second integral on the right hand side can be computed by using equations \eqref{eq:maincalcg1} and \eqref{eq:hedegen},
\be
 \int_{Y}\kappa_1\psi^{2g-2}=\int_{\Mb_{1,1}}\kappa_1 \cdot \int_{\Hb_{g-1,1}^1}\psi^{2g-2} 
 =\frac{1}{24}\cdot \frac{1}{2^{2g-2}(2g-1)!}.
\ee

Adding this to the value of the integral \eqref{eq:maincalc}, we obtain
\begin{multline}
\int_{\Mb_{g,1}} \frac{ 1 - \lambda + \lambda^2 - \dotsm  \pm \lambda^{g-1}}{ (1-\psi)(1-2\psi)}  \kappa_1\psi_1^{2g-2} \\ =
\frac{12g^2 -12g +3 +2g^2 +g}{3\cdot2^{2g}(2g+1)!}  =  \frac{14g^2 -11g +3  }{3\cdot{}2^{2g}(2g+1)!}.
\end{multline}
which is equivalent to the statement of the corollary.

Finally, it is interesting to point out that by equation \eqref{eq:maincalcg1}, the corollary also holds for $g=1$ if the integral is taken in the geometric, rather than orbifold, sense.  

\section{Concluding Remarks}

As mentioned in the introduction, the main calculation of this paper is the first direct  application of Witten's duality, and in some sense the Witten cycles themselves,  to questions about the tautological ring.  The fact that such an application has not appeared earlier does not necessarily put into question  the significance of  the Witten cycles.  On the contrary,  the intimate connections between the cycles $W_a$ and  matrix models have  been known (or at least suspected)  since the introduction of the combinatorial cycles and seem to reflect a deep interplay between mathematics and physics.  

 The lack of an earlier application of this perspective does, however, emphasize the difficulties in using chain level intersection techniques to study the moduli space of curves, difficulties due primarily  to the highly degenerate nature of the intersections that  usually accompany combinatorial cycles.
 In fact, the calculation of this paper is a bit of a lucky coincidence where the chain level intersections of cycles, while not being generic, are still non-degenerate.  


Adding to  the lucky coincidence, the case of Witten's duality  that we utilize here is particularly simple.  Indeed, the $a=1$ case given by  \eqref{eq:introduality} exhibits a 
 clear distinction between the combinatorics of the open moduli (given by $W_1$) and that of the boundary $\partial\M_{g,1}$.  For $a>1$, this distinction is not so clear.  Indeed, complete descriptions of the boundary classes $\delta_a$ in theorem \ref{thm:generalduality} have still not been obtained for all $a$, and Arbarello and Cornalba \cite{arbarello} have shown that any such description will likely exhibit some essential ambiguities.

Still, the classes $\delta_a$ can be described for small $a$, and it seems reasonable to try to use our techniques to attack the integrals\footnote{ We note that 
 j. Bertin and M. Romagny have  determined these integrals \eqref{eq:generalizations} by algebraic means \cite{bertin}.  }
\be\label{eq:generalizations}
\int_{\Hgb^1} \kappa_a\psi^{2g-1-a}, 
\ee
if not for general $a$, at least for these few cases. 
The main obstacle to computing these integrals  is that for $a>1$, the intersection $\Hgc\cap W_a$  is no longer non-degenerate.  This can be seen by noticing that there exist (real) codimension $k$ cells in $\Hgc$ which contain metric fatgraphs with a $2k+2$-valent vertex fixed under the hyperelliptic involution.  The case with $k=2$ gives the non-degenerate intersection $\Hgc\cap W_1$;  however, for $k>2$ the intersection is obviously degenerate. 

 Regardless, it may still be feasible to handle degenerate intersections by a more careful analysis of the cell structure of the higher dimensional Stasheff polytopes, and if so, it could serve as a positive step towards handling the more difficult degenerate intersections of the Witten cycles amongst themselves. 

As a final note, there is also an algebro-geometric obstacle  to obtaining relations among Hodge integrals from the integrals \eqref{eq:generalizations} for $a>2$, as the higher dimensional components of the degeneracy locus in Porteus' formula no longer give vanishing contributions in these cases.

\emph{Acknowledgments.} 
A special thanks is owed to the author's advisor Kefeng Liu for first introducing him to the moduli space of curves.  The author is also 
deeply indebted to Robert Penner for his enthusiasm and support throughout all stages of this research.  Conversations with Geoffrey Mess and Ravi Vakil have also been very helpful in the preparation of this paper.  Finally,  the author is appreciative of the hospitality shown by  the Centre for Mathematical Sciences in  Hangzhou, China  during the summer of 2004 where the main relations of this paper were first derived.

\small
\textsc{Department of Mathematics, University of Southern California,}

\textsc{Los Angeles, California 90089-2532}

\emph{E-mail address}: bene@usc.edu

\end{document}